\documentclass[envcountsame]{llncs}
\usepackage{amsmath,amsfonts,amssymb}
\usepackage{algpseudocode}
\usepackage[all]{xy}

\newcommand{\CC}{\bbbc}
\newcommand{\FF}{\bbbf}
\newcommand{\MM}{\bbbm}

\newcommand{\QQ}{\bbbq} 
\newcommand{\ZZ}{\bbbz}

\newcommand{\OO}{\mathcal{O}}

\newcommand{\dsp}{\displaystyle}
\newcommand{\Div}{\mathrm{Div}}
\newcommand{\End}{\mathrm{End}}
\newcommand{\Jac}{\mathrm{Jac}}

\newcommand{\subgrp}[1]{\langle{#1}\rangle}

\newcommand{\isom}{\cong}

\newcommand{\genus}[1]{{g_{#1}}}
\newcommand{\variety}[1]{{V({#1})}}
\newcommand{\closure}[1]{{\overline{#1}}}

\spnewtheorem{algorithm}[theorem]{Algorithm}{\bfseries}{\rm}

\begin{document}

\title{Efficiently Computable Endomorphisms\\ for Hyperelliptic Curves}
\author{David R.~Kohel, Benjamin A.~Smith}
\authorrunning{D.~R.~Kohel, B.~A.~Smith}
\institute{
	School of Mathematics and Statistics, \\
	The University of Sydney \\
	\email{\{kohel,bens\}@maths.usyd.edu.au}
}

\maketitle


\begin{abstract}
Elliptic curves have a well-known and explicit theory for the 
construction and application of endomorphisms, which can be applied 
to improve performance in scalar multiplication.  Recent work has 
extended these techniques to hyperelliptic Jacobians, but one 
obstruction is the lack of explicit models of curves together 
with an efficiently computable endomorphism.  In the case of 
hyperelliptic curves there are limited examples, most methods 
focusing on special CM curves or curves defined over a small 
field.  In this article we describe three infinite families of 
curves which admit an efficiently computable endomorphism, and 
give algorithms for their efficient application.  

\end{abstract}

\section{Introduction}

The use of efficiently computable endomorphisms for speeding up 
point multiplication on elliptic curves is well-established for 
elliptic curves and more recently has been used for hyperelliptic 
curves.  Koblitz~\cite{Koblitz-1991} proposed $\tau$-adic expansions 
of the Frobenius endomorphism on curves over a small finite fields.
Gallant, Lambert, and Vanstone~\cite{GLV-2001} later proposed 
using an expression 
$$
	[k]P = [k_0]P + [k_1]\phi(P)
$$
on more general curves to evaluate multiplication by $k$ on 
a point $P$, using an efficiently computable endomorphism $\phi$. 
Various improvements and combinations of these methods have been 
proposed for both elliptic and hyperelliptic 
curves~\cite{Lange-2001,Solinas-2000,CLSQ-2003}.

One feature of elliptic curves, not available for multiplicative 
groups of finite fields, is the freedom to choose a parameter:
geometrically they form a one-dimensional family,
parametrized by the $j$-invariant.  
Restriction to curves of a special form 
destroys this degree of freedom.  While no proof exists that 
special curves, CM curves or Koblitz curves are less insecure, 
these nonrandom curves can be qualitatively distinguished from 
their nonrandom cousins in terms of their endomorphism rings. 
Thus preference is often given to curves randomly selected over 
a large finite field when performance is not the determining 
issue. 

In contrast, hyperelliptic curves of genus $g$ admit a much 
larger degree of freedom.  
In genus $2$, they form a three dimensional family:
curves with different classifying triple of invariants $(j_1,j_2,j_3)$ 
can not be isomorphic over any extension field.  
Until the recent work of Takashima~\cite{Takashima-2006}, 
the only curves proposed for 
cryptographic use with efficiently computable endomorphisms 
are either the CM curves with exceptional automorphisms --- 
the analogues of elliptic curves $y^2 = x^3 + a$ or $y^2 = 
x^3 + ax$ --- or Koblitz curves --- curves defined over a small 
field with point on the Jacobian taken over a large prime 
degree extension (see Park et al.~\cite{PSL-2002} for the former 
and Lange~\cite{Lange-2001} for the latter).
Besides the notable exceptions of CM curves with exceptional 
automorphisms, curves with CM have been exploited for point 
counting but not for their endomorphism ring structure, for 
lack of a constructive theory of efficiently computable 
endomorphisms.  

In this work,
we address the problem of effective algorithms for 
endomorphisms available on special families of curves.  
We describe three families, of dimensions $1$, $1$, and $2$
respectively, 
of curves whose Jacobians admit certain {\it real} endomorphisms.  
First, we introduce the general framework for 
constructing endomorphisms via correspondences derived from 
covering curves.  
Subsequently, 
we provide a one-dimensional family 
derived from Artin--Schreier covers, 
then describe a construction of Tautz, Top, and Verberkmoes~\cite{TTV-1991} 
for a one-dimensional family of curves with explicit 
endomorphisms deriving from cyclotomic covers.  
Finally,
we describe an elegant construction of Mestre~\cite{Mestre-1991} 
from which we obtain a two-dimensional family of curves whose 
Jacobians admit explicit endomorphisms, derived from covers 
of elliptic curves.  
In each case we develop explicit 
algorithms for efficient application of the endomorphism, 
suitable for use in a GLV decomposition.  
Independently, Takashima~\cite{Takashima-2006} provided an 
efficient algorithm for endomorphisms in the latter family 
(in terms of variants of Brumer and Hashimoto) with real 
multiplication by $(1+\sqrt{5})/2$.  
These families provide
a means of generating curves randomly selected within a large family,
yet which admit efficiently computable endomorphisms.  
 
\section{Arithmetic on Hyperelliptic Jacobians}

In the sequel we denote by $X/k$ 
a hyperelliptic curve of genus $\genus{X}$
in the form
$$
	v^2 = f(u) = u^{2\genus{X}+1} + c_{2\genus{X}}u^{2\genus{X}} + \cdots +
c_0,$$
with each $c_i$ in $k$,
which we require to be a field of characteristic not $2$.
The Jacobian of $X$, denoted $\Jac(X)$, 
is a $\genus{X}$-dimensional variety 
whose points form an abelian group.
Let $\OO$ denote the point at infinity of $X$.
Each point $P$ on $\Jac(X)$
may be represented by a 
divisor on $X$, 
that is, as a formal sum of points 
$$
	P 
	= \sum_{i=1}^{m} [P_i] - m [\OO]
  	= \sum_{i=1}^{m} [(u_i,v_i)] - m [\OO],
$$
where $m \le \genus{X}$.  
We say such a divisor is {\em semi-reduced}
if $(u_i,v_i)\not=(u_j,-v_j)$ for all $i\not= j$.
For a point to be defined in $\Jac(X)(k)$,
its divisor must be Galois-stable;
the representation as a divisor
has the disadvantage that the individual points $(u_i,v_i)$
may be defined only over some finite extension $K/k$.
Thus, for computations, 
we use instead the Mumford representation for divisors, 
identifying $P$ with the 
ideal class
$$
	P = [(a(u),v-b(u))],
$$
where $a$ and $b$ are polynomials in $k[u]$
such that 
$ a(u) = \prod_i(u-u_i) $ 
and $ v_i = b(u_i) $ for all $i$.  
In this guise, 
addition of points $P$ and $Q$ is an ideal product, 
followed by a reduction algorithm to 
produce a unique ``reduced'' ideal representing $P+Q$.  
Cantor~\cite{Cantor-1987} provides algorithms 
to carry out these operations.

\begin{algorithm}
\label{CantorReduction}
	Given a semi-reduced representative $(a(u),v - b(u))$
	for a point $P$ on the Jacobian of a hyperelliptic curve
	$X: v^2 = f(u)$, 
	returns the reduced representative of $P$.
\begin{algorithmic}
\Function{CantorReduction}{$(a(u), v - b(u))$}
	\While{$\deg(a) > \genus{X}$}
	\State $a := (f-b^2) / a$;
	\State $b := -b \bmod a$;
	\EndWhile;
	\State $a := a/\textsc{LeadingCoefficient}(a)$;
	\State \textbf{return} $(a, v - b(u))$;
\EndFunction;
\end{algorithmic}
\end{algorithm}

\noindent
Each iteration of Algorithm~\ref{CantorReduction}
replaces $a$ with a polynomial of degree
$\mathrm{max}(2\genus{X}+1-\deg(a), \deg(a) - 2)$.
It follows that Algorithm~\ref{CantorReduction}
will produce a reduced representative 
for the ideal class $[(a(u),v-b(u))]$
after $\lceil (\deg(a) - \genus{X})/2 \rceil$ iterations.

\section{Explicit Endomorphisms}

Let $C$ be a curve with an automorphism $\zeta$,
and let $\pi: C \rightarrow X$ 
be a covering of $X$.  
We have two coverings,
$\pi$ and $\pi\circ\zeta$,
from $C$ to $X$;
together, they
induce a map $\eta$ of divisors 
$$
	\eta := (\pi\circ\zeta)_*\pi^*: \Div(X) \rightarrow \Div(X),
$$
where 
$$
\pi^*([P]) = \!\!\!\sum_{Q \in \pi^{-1}(P)}\!\!\! e_\pi(Q)[Q]
\mbox{\ \ and\ \ }
(\pi\circ\zeta)_*([Q]) = [\pi(\zeta(Q))].
$$
This map on divisors induces an endomorphism of the Jacobian $\Jac(X)$,
which we also denote $\eta$.

In our constructions,
we take $\pi$ to be the quotient 
by an involution $\sigma$ of $C$,
so that $\pi$ is a degree-$2$ covering,
and $\pi = \pi\circ\sigma$.
Thus 
$$
	\pi^*([P]) = [Q] + [\sigma(Q)]
$$
for any point $Q$ in $\pi^{-1}(P)$.  
We will take
$\zeta$ to be an automorphism of $C$ of prime order $p$,
such that $\subgrp{\zeta,\sigma}$
is a dihedral subgroup
of the automorphism group of $C$:
that is, $\sigma\zeta = \zeta^{-1}\sigma$.
The following proposition describes
the resulting endomorphism 
$\eta = (\pi\circ\zeta)_*\circ\pi^*$.

\begin{proposition}
\label{dihedralendomorphisms}
	Let $C$ be a curve
	with an involution $\sigma$
	and an automorphism $\zeta$
	of prime order $p$
	such that $\sigma\zeta = \zeta^{-1}\sigma$.
	Let $\pi: C \to X := C/\subgrp{\sigma}$
	be the quotient of $C$ by the action of $\sigma$,
	and let $\eta := (\pi\circ\zeta)_*\circ\pi^*$
	be the endomorphism of $\Jac(X)$ induced by $\zeta$.
	The subring $\ZZ[\eta]$ of $\End(\Jac(X))$
	is isomorphic to $\ZZ[\zeta_p + \zeta_p^{-1}]$,
	where $\zeta_p$ is a primitive $p^\mathrm{th}$
	root of unity over $\QQ$.
\end{proposition}
\begin{proof}
	The subring $\ZZ[\zeta_* + \zeta_*^{-1}]$ of $\Jac(C)$
	is isomorphic to $\ZZ[\zeta_p + \zeta_p^{-1}]$,
	since $p$ is prime.
	The statement follows upon noting that
	the following diagram commutes.
	\SelectTips{cm}{10}
	$$
	\xymatrix@R=16mm@C=20mm{
	  \Jac(C) \ar[r]^{\dsp\zeta_* + \zeta^{-1}_*} 
                  \ar[d]_{\dsp\pi_*} & \Jac(C) \ar[d]^{\dsp\pi_*} \\
	  \Jac(X) \ar[r]^{\dsp\eta = \pi_*\zeta_*\pi^*} & \Jac(X)
	}
	$$
	To see this, 
	observe that for any $Q$ in $\Jac(C)$ we have
	$$ \begin{array}{r@{\;=\;}l}
		\eta(\pi_*(Q)) & \pi_*\zeta_*\pi^*\pi_*(Q) \\
				& \pi_*\zeta_*(1 + \sigma_*)(Q) \\
				& \pi_*(\zeta_* + \sigma_*\zeta^{-1}_*)(Q) \\
				& \pi_*(\zeta_* + \zeta^{-1}_*)(Q) , \\
	\end{array}
	$$
	since $\pi^*\pi_* = (1 + \sigma_*)$ and $\pi_*\sigma_* = \pi_*$.
	See also Ellenberg \cite[\S2]{Ellenberg}.
\qed
\end{proof}

\noindent
Suppose $C$, $X$, $\pi$, $\zeta$ and $\eta$
are as in Proposition~\ref{dihedralendomorphisms}.
Our aim is to give an explicit realization of 
the endomorphism $\eta$ of $\Jac(X)$,
in the form of a map on ideal classes.
To do this,
we form the algebraic correspondence
$$
	Z := (\pi\times(\pi\circ\zeta))(C) \subset X\times X.
$$
Let $\pi_1$ and $\pi_2$ be
the restrictions to $Z$
of the projections from $X\times X$
to its first and second factors, respectively;
then
$\eta = (\pi_2)_*\circ\pi_1^*$.
We will give an affine model for $Z$
as the variety cut out by an ideal in 
$k[u_1,v_1,u_2,v_2]/(v_1^2 - f(u_1),v_2^2 - f(u_2))$;
for this model,
the maps $\pi_1$ and $\pi_2$
are defined by $\pi_i(u_1,v_1,u_2,v_2) = (u_i,v_i)$.

Suppose that $Z$ is defined by an ideal $(v_2 - v_1, E(u_1,u_2))$,
where $E$ is quadratic in $u_1$ and $u_2$
(this will be the case in each of our constructions).
If $(u,v)$ is a generic point on $X$, then $\pi_1^*([(u,v)])$ 
is the effective divisor on $Z$ cut out by $(v_2 - v, E(u,u_2))$.
Therefore, if $e_1$ and $e_2$ are the solutions in $\closure{k(u)}$
to the quadratic equation $E(u,x) = 0$ in $x$,
then 
$$
	\eta([(u,v)]) 
	= (\pi_2)_*\pi_1^*([(u,v)]) 
	= [(e_1, v)] + [(e_2,v)].
$$
It remains to translate this description of the action of $\eta$
in terms of points into a map on ideal classes.

Suppose $[(a(u),v-b(u))]$ is a point on $\Jac(X)$.
Extending the above, we have
$$
\begin{array}{r@{\;=\;}l}
	\eta([(a(u),v-b(u))])
  	& [(a(e_1),v-b(e_1))] + [(a(e_2),v-b(e_2))] \\
  	& \displaystyle [(N(a), v-\frac{(f(u)+N(b))}{T(b)} \bmod N(a)],
\end{array}
$$
where $N(a) = a(e_1)a(e_2)$, 
$N(b) = b(e_1)b(e_2)$, and 
$T(b) = b(e_1) + b(e_2)$.\footnote{%
	The modular inversion of $T(b)$ should be carried out 
	after clearing denominators
	and removing common factors from $N(a)$, $T(b)$, and 
	$f(u)+N(b)$ (generically, $N(a)$ and $T(b)$ are coprime).
	Proposition~\ref{etaimageproposition} below makes this precise.
}  
Since functions $T(a)$, $N(b)$ and $T(b)$ 
are symmetric polynomials in $e_1$ and $e_2$, 
we can write each as a polynomial in 
the elementary symmetric functions $e_1 + e_2$ and $e_1e_2$.
Moreover, 
$e_1 + e_2$ and $e_1e_2$
are elements of $k(u)$:
if $E(u,x) = E_2(u)x^2 + E_1(u)x + E_0(u)$, 
then 
$e_1 + e_2 = -E_1/E_2$ and $e_1e_2 = E_0/E_2$.

\begin{definition}
\label{TNdefns}
	For any polynomial $a(x)$ over $k$, 
	we define 
	$T(a) = a(e_1)+a(e_2)$ and 
	$N(a) = a(e_1)a(e_2)$, 
	and for $i,j \ge 0$ we define 
	$$
		t_i := e_1^i + e_2^i,\ n_i := (e_1e_2)^i
		\quad \mbox{and}\quad 
		n_{i,j} := e_1^ie_2^j + e_1^je_2^i.
	$$
\end{definition}

\noindent
Note that $t_i$ and $n_{ij}$ are elements of $k(u)$ and that
\begin{equation}
\label{Texpression}
	T\left(\sum_{i=0}^\genus{X} a_i x^i\right) 
	= \sum_{i=0}^\genus{X} a_i t_i,
\end{equation}
and
\begin{equation}
\label{Nexpression}
	N\left(\sum_{i=0}^\genus{X} a_i x^i\right) 
	= \sum_{i=0}^\genus{X}\sum_{i=0}^\genus{X} a_ia_j n_{i,j}.
\end{equation}
The following elementary lemma provides simple recurrences for 
the construction of the sequences $\{t_i\}$ and $\{n_{i,j}\}$.

\begin{lemma}
\label{TNrecur}
The elements $t_i$, $n_i$ and $n_{i,j}$ 
satisfy the following recurrences:
\begin{enumerate}
	\item 	$n_{i+1} = (e_1e_2)n_i$ for $i \ge 0$, 
		with $n_0 := 1$;
	\item 	$t_{i+1} = (e_1+e_2)t_i - (e_1e_2)t_{i-1}$ 
		for $i \ge 1$, 
		with $t_0 = 2$ and $t_1 = (e_1 + e_2)$;
	\item 	$n_{i,i} = n_i$ and $n_{i,j} = n_it_{j-i}$ 
		for $i \ge 0$ and $j > i$.
\end{enumerate}
\end{lemma}

\noindent
Equations~\eqref{Texpression} and~\eqref{Nexpression} above
express $T$ and $N$ in terms of the functions $t_i$ and $n_{i,j}$,
which depend \emph{only} upon $t_1$ and $n_1$
by Lemma~\ref{TNrecur}.
Thus,
given 
$t_1 = e_1 + e_2$ and $n_1 = n_{1,1} = e_1e_2$, 
the recurrences of Lemma~\ref{TNrecur}
give a simple and fast algorithm 
for computing the maps $T$ and $N$.
If we further assume that $T$ and $N$ will only be evaluated
at polynomials $a$ and $b$ from reduced ideal class representatives
$(a(u),v-b(u))$,
then we need only compute the $t_i$ and $n_{i,j}$ 
for $0 \le i \le j \le \genus{X}$.

\begin{algorithm}
Given functions $t_1$ and $n_1$ in $k(u)$, 
together with the genus $\genus{X}$ of a curve $X$, 
returns the maps $T$ and $N$ of Definition \ref{TNdefns}.
\label{TNalgo}
\begin{algorithmic}
\Function{RationalMaps}{$t_1$,$n_1$,$\genus{X}$}
	\State 	$n_0 := 1$;
	\State  $t_0 := 2$;
	\For{$i$ \textbf{in} $[1,\ldots,\genus{X}]$}
		\State $n_{i+1} := n_1 n_i$;
		\State $t_{i+1} := t_1t_i - n_1t_{i-1}$;
	\EndFor;
	\For{$i$ \textbf{in} $[1,\ldots,\genus{X}]$}
		\State $n_{i,i} := n_i$;
		\For{$j$ \textbf{in} $[i+1,\ldots,\genus{X}]$}
			\State $n_{i,j} := n_it_{j-i}$;
		\EndFor;
	\EndFor;
	\State $T := (\sum_{i=0}^\genus{X}a_iX^i 
			\longmapsto 
			\sum_{i=0}^\genus{X}a_it_i)$;
	\State $N := (\sum_{i=0}^\genus{X}a_iX^i 
			\longmapsto 
			\sum_{i=0}^\genus{X}\sum_{j=i}^\genus{X}a_ia_jn_{i,j})$;
	\State \textbf{return} $T$, $N$;
\EndFunction;
\end{algorithmic}
\end{algorithm}

\noindent
The following proposition shows that the maps $T$ and $N$
may be used 
to compute $\eta([(a(u),v-b(u))])$ for all points $[(a(u),v-b(u))]$ 
of $\Jac(X)$.

\begin{proposition}
\label{etaimageproposition}
	Let $\eta$ be the endomorphism of $\Jac(X)$
	induced by a correspondence
	$\variety{v_2 - v_1, E_2(u_1)u_2^2 + E_1(u_1)u_2 + E_0(u_1)}$
	on $X\times X$;
	set $t_1 = -E_1/E_2$ and $n_1 = E_0/E_2$,
	and let $T$ and $N$ be the maps of Definition~\ref{TNdefns}.
	If $(a(u),v-b(u))$ is the reduced representative 
	of a point $P$ of $\Jac(X)$,
	then $\eta(P)$
	is represented by
	$$
		\left(
			\frac{E_2^\genus{X} N(a)}{G},
			v - \left( \frac{(f + N(b))/G}{T(b)/G} 
				\bmod \frac{E_2^\genus{X}N(a)}{G} 
			\right)
		\right),
	$$
	where $G = \gcd(E_2^\genus{X}N(a),E_2^\genus{X}T(b))$.
	Algorithm~\ref{CantorReduction}
	computes the reduced representative of $\eta(P)$
	after at most $\lceil\genus{X}/2\rceil$ iterations 
	of its main loop.
\end{proposition}
\begin{proof}
	We have 
	$$
	\begin{array}{r@{\;=\;}l}
		\eta([(a(u),v-b(u))]) 
		& [(a(e_1),v-b(e_1))(a(e_2),v-b(e_2))] \\
		& [(N(a), v^2 - T(b)v + N(b))] \\
		& [(E_2^\genus{X})(N(a), T(b)v - (f + N(b)))]. \\
		& [(E_2^\genus{X}N(a), E_2^\genus{X}T(b)v - E_2^\genus{X}(f +
N(b)))]. \\	\end{array}
	$$
	It is easily verified that
	$E_2^\genus{X}N(a)$, 
	$E_2^\genus{X}T(b)$ 
	and $E_2^\genus{X}(f + N(b))$ are polynomials,
	and that if $G = \gcd(E_2^\genus{X}N(a),E_2^\genus{X}T(b))$,
	then $G$ also divides $E_2^\genus{X}(f + N(b))$.
	Therefore
	$$
	\begin{array}{r@{\;=\;}l}
		\eta([(a(u),v-b(u))]) 
		& [(G)(E_2^\genus{X}N(a)/G, E_2^\genus{X}T(b)v/G -
E_2^\genus{X}(f + N(b))/G))] \\
		& [(E_2^\genus{X}N(a)/G,
E_2^\genus{X}T(b)v/G - E_2^\genus{X}(f + N(b))/G)] \\		&
[(E_2^\genus{X}N(a)/G, v - I\cdot E_2^\genus{X}(f + N(b))/G)], \\
\end{array}	$$
	where $I$ denotes 
	the inverse of $E_2^\genus{X}(f + N(b))/G$ modulo
$E_2^\genus{X}N(a)/G$,	proving the first claim.
	Now, 
	if $(a(u), v - b(u))$ is the reduced representative of $P$,
	then $\deg(a) \le \genus{X}$,
	so the degree of $E_2^\genus{X}N(a)$ is at most $2\genus{X}$.
	After each iteration of Algorithm~\ref{CantorReduction},
	the degree of $a$ becomes 
	$\mathrm{max}(2\genus{X}+1-\deg(a), \deg(a) - 2)$,
	and the algorithm terminates
	when $\deg(a) \le \genus{X}$;
	this occurs after $\lceil\genus{X}/2\rceil$ iterations.
\qed
\end{proof}

\noindent
The following algorithm applies Proposition~\ref{etaimageproposition}
to compute the image of a point of $\Jac(X)$ under $\eta$.
This gives an explicit 
realization of $\eta$ as a map on ideal classes.

\begin{algorithm}
\label{evaluationalgo}
	Given a point $P$ on the Jacobian of a curve $X: v^2 = f(u)$
	and rational maps $T$ and $N$ derived for an endomorphism
	$\eta$ of $\Jac(X)$ using Algorithm~\ref{TNalgo},
	returns the reduced ideal class representative of $\eta(P)$.
\begin{algorithmic}
\Function{Evaluate}{$P = (a(u), v-b(u))$, $T$, $N$}
	\State $a' := N(a)$;
	\State $d := T(b)$;
	\State $E :=
\textsc{LCM}(\textsc{Denominator}(a'),\textsc{Denominator}(d))$;
	\State
$G := \textsc{GCD}(\textsc{Numerator}(a'),\textsc{Numerator}(d))$;
	\State $a' := E\cdot a'/G$;
	\State $d := E \cdot d/G$;
	\State $I := d^{-1} \pmod{a'}$;
	\State $b' := I\cdot E\cdot(f + N(b))/G \pmod{a'}$;
	\State \textbf{return} \textsc{CantorReduction}($(a', v - b')$);
\EndFunction;
\end{algorithmic}
\end{algorithm}

\begin{remark}
	In the families of curves described below in 
	Sections~\ref{artinschreiersection} 
	and~\ref{cyclotomicsection}
	below, 
	$T$ and $N$ are polynomial maps,
	and we may take $E = 1$ in Algorithm~\ref{evaluationalgo}.
\end{remark}

\section{Applications I: Curves with Artin--Schreier Covering}
\label{artinschreiersection}

In this section 
we construct a family of curves $X_p$ 
in one free parameter $t$ 
for each prime $p \ge 5$, 
and determine explicit 
endomorphisms deriving from a cover by the Artin--Schreier curve 
defined over $\FF_p$ by  
$$
	C_p : y^p - y = x + \frac{t}{x}\cdot
$$
The eigenvalues of Frobenius in this family are described by 
classical Kloosterman sums~\cite{Weil-1948}.  

An analogous family $y^2 = x^p - x + t$ was described by Duursmaa 
and Sakurai~\cite{DuuSak00}, for which the automorphism 
$x \mapsto x+1$ was proposed for efficient scalar multiplication.  
In constrast to our family, every member of this family is isomorphic 
over a base extension to the supersingular curve $y^2 = x^p - x$.

\subsection{Construction of the Artin--Schreier Covering}

The curve $C_p$ has automorphisms $\zeta$ (of order $p$) and 
$\sigma$ (of order $2$), defined by
$$
	\zeta(x,y) = (x, y+1) 
	\quad \mbox{and} \quad
	\sigma(x,y) = (-t/x,-y).
$$
Let $X_p$ be the quotient of $C_p$ by $\subgrp{\sigma}$, with 
affine model
$$
	X_p: v^2 = f(u) = u(u^{(p-1)/2} - 1)^2 - 4t.
$$
The quotient map $\pi: C_p \rightarrow X_p$ is a covering of 
degree $2$, sending $(x,y)$ to $(u,v) = (y^2,x-t/x)$.
Observe that $X_p$ is a family of curves of genus $(p-1)/2$.

The automorphism $\zeta$ of $C_p$ induces an endomorphism
$\eta := (\pi\circ\zeta)_*\pi^*$ on $\Jac(X_p)$, whose minimal 
polynomial equals that of $\eta_p = \zeta_p+\zeta_p^{-1} \in \CC$.  
The endomorphism $\eta$ is induced by the correspondence 
$Z := (\pi\circ\zeta\times\pi)(C_p)$ on $X_p\times X_p$, for 
which we may directly compute an affine model
$$
	Z = \variety{
		v_2 - v_1,
		u_2^2 + u_1^2 - 2u_1u_2 - 2u_2 - 2u_1 + 1
	}.
$$
Setting $t_1 := 2(u+1)$ and $n_1 := (u - 1)^2$ and applying 
Algorithm~\ref{TNalgo}, 
we obtain polynomial maps $T$ and $N$ such that
$\eta$ is realized by $P \mapsto \textsc{Evaluate}(P,T,N)$, 
using Algorithm~\ref{evaluationalgo}.
The first few $t_i$ and $n_{i,j}$ derived in Algorithm~\ref{TNalgo}
are given in Table~\ref{artinschreiertsandns} below.

\begin{proposition}
\label{prop_endo_artin-schreier}
The Jacobian $\Jac(C_p)$ is isogenous to $\Jac(X_p)^2$, and its 
endomorphism ring contains an order in $\MM_2(\QQ(\eta_p))$.
\end{proposition}

\begin{proof}
The automorphisms $\zeta$ and $\sigma$ determine a homomorphic image of the 
group algebra $A = \QQ[\langle\zeta,\sigma\rangle]$ in $\End^\circ(\Jac(C_p))$.
But $A$ is a semisimple algebra of dimension $2p$ over $\QQ$, whose simple 
quotients are of dimensions $1$, $1$, and $2\varphi(p)$.  
Moreover,
$\zeta + \zeta^{-1}$ is in the centre of $A$ 
and generates a subring isomorphic to $\QQ \times \QQ(\eta_p)$.  
Since $\zeta$ and $\sigma$ do not commute, it follows that the latter 
algebra is isomorphic to $\MM_2(\QQ(\eta_p))$.

Let $e_1$ and $e_2$ be the central idempotents associated to the quotients 
of dimensions $1$.  On each associated abelian variety $e_i\Jac(C_p)$, the 
automorphism $\zeta$ acts trivially, thus maps through the Jacobian of the 
genus~$0$ quotient $C_p/\langle\zeta\rangle$; it follows that the image of 
$A$ in $\End^\circ(\Jac(C_p))$ is isomorphic to $\MM_2(\QQ(\eta_p))$.

Let $\epsilon_1 = 1+\sigma$ and $\epsilon_2 = 1-\sigma$.  Noting that 
$$
\epsilon_i^2 = 2\epsilon_i,\quad 
\epsilon_1\epsilon_2 = 0, \mbox{ and }
\epsilon_1 + \epsilon_2 = 2,
$$
we let $A_1 = {\epsilon_1}_*\Jac(C_p)$ and $A_2 = {\epsilon_2}_*\Jac(C_p)$ 
be subabelian varieties of $\Jac(C_p)$ such that $\Jac(C_p) = A_1 + A_2$, 
and $A_1 \cap A_2$ is finite.  Since $\zeta-\zeta^{-1}$ determines an 
isogeny $\psi = \zeta_* - \zeta^{-1}_*$ of $\Jac(C_p)$ to itself, the 
relation
$$
(\zeta-\zeta^{-1})\epsilon_1 = \epsilon_2(\zeta-\zeta^{-1}),
$$
implies that 
$
\psi(A_1) = {\epsilon_2}_* \psi(\Jac(C_p)) = A_2,
$
so that $A_1$ and $A_2$ are isogenous.  
But $\pi_*$ is an isogeny of $A_1$ to $\Jac(X_p)$, 
whence $\Jac(C_p) \sim \Jac(X_p)^2$. 
\qed
\end{proof}

\begin{corollary}
The Jacobian $\Jac(X_p)$ has a rational $p$-torsion point.
In particular,
$\Jac(X_p)$ is not a supersingular abelian variety.
\end{corollary}

\begin{proof}
The curve $C_p$ has two rational points fixed by $\zeta$, whose 
difference determines a point in $\ker(1-\zeta_*)$.  
But 
$$
(1-\zeta)(1-\zeta^2)\cdots(1-\zeta^{p-1}) = p,
$$
so $\ker(1-\zeta_*)$ is contained in $\Jac(X_p)[p]$.
If $\chi(T)$ and $\xi(T)$ are the characteristic polynomials of 
Frobenius on $\Jac(C_p)$ and $\Jac(X_p)$, respectively, then 
$\chi(T) = \xi(T)^2$.  Since $|\Jac(C_p)(k)| = \chi(1)$ is 
divisible by $p$, so is $|\Jac(X_p)(k)| = \xi(1)$.
\qed
\end{proof}

\begin{remark}
In fact, it is possible to show that the $p$-rank of $\Jac(X_p)$ 
is exactly equal to $1$,
so the Jacobians are neither ordinary 
nor supersingular.
\end{remark}

\begin{table}
\caption{Artin--Schreier covers: $t_i$ and $n_{i,j}$ for $0 \le i \le j
\le 3$.}\label{artinschreiertsandns}
\vspace{-4mm}
\begin{center}
\begin{tabular}{|c|l|}
\hline
$t_0$		& $2$	\\
$t_1$		& $2(u+1)$	\\
$t_2$		& $2(u^2 + 6u + 1)$	\\
$t_3$		& $2(u^3 + 15u^2 + 15u + 1)$	\\
\hline
$n_{0,0}$ 	& $1$ \\
$n_{0,1}$ 	& $2(u+1)$ \\
$n_{0,2}$ 	& $2(u^2 + 6u + 1)$ \\
\hline
\end{tabular}
\begin{tabular}{|c|@{\;}l|}
\hline
$n_{0,3}$	& $2(u^3 + 15u^2 +15u + 1)$ \\
$n_{1,1}$ 	& $(u-1)^2$ \\
$n_{1,2}$ 	& $2(u-1)^2(u+1)$ \\
$n_{1,3}$	& $2(u-1)^2(u^2 + 6u + 1)$ \\
$n_{2,2}$ 	& $(u-1)^4$ \\
$n_{2,3}$	& $2(u-1)^4(u+1)$ \\
$n_{3,3}$	& $(u-1)^6$ \\
\hline
\end{tabular}
\end{center}
\vspace{-4mm}
\end{table}

\subsection{Hyperelliptic Curves of Genus $2$ with Real Multiplication
by $\eta_5$}

For $p = 5$, the construction above yields a one-parameter family of
genus $2$ hyperelliptic curves defined by
$$
	X_5: v^2 = f_5(u) = u(u^2-1)^2 + t,
$$
whose Jacobian has endomorphism ring containing 
$\ZZ[\eta_5] \isom \ZZ[x]/(x^2+x-1)$.

Each point $P$ of $\Jac(X_5)$ 
may be represented by an ideal $(a(u),v-b(u))$
with $a$ and $b$ of degrees $2$ and $1$ respectively:
hence, suppose
$a(u) = a_2 u^2 + a_1 u + a_0$ 
and $b(u) = b_1 u + b_0$. 
Applying Algorithm~\ref{TNalgo},
we see that
$$
\begin{array}{r@{\;=\;}l}
	N(a) & a_2^2n_{2,2} + a_2a_1n_{1,2} + a_1^2n_{1,1}
         	+ a_2a_0n_{0,2} + a_1a_0n_{0,1} + a_0^2n_{0,0}, \\
	N(b) & b_1^2 n_{1,1} + b_1b_0 n_{0,1} + b_0^2 n_{0,0},
	\; \mbox{and} \\
	T(b) & 2b_1(u+1) + 2b_0 ,
\end{array}
$$
with the $n_{i,j}$ as in Table~\ref{artinschreiertsandns}.
The endomorphism $\eta$ is then explicitly realized by
$\eta: P \mapsto \textsc{Evaluate}(P,T,N)$,
using Algorithm~\ref{evaluationalgo}.

\begin{remark}
	The Igusa invariants of the curve $X_5$ 
	determine the weighted projective point 
	$(J_2:J_4:J_6:J_8:J_{10}) = (3:2:0:4:4t^2)$. 
	In particular,
	the curves determine a one-dimensional subvariety
	of the moduli space of genus $2$ curves.
\end{remark}

\subsection{Hyperelliptic Curves of Genus $3$ with Real Multiplication
by $\eta_7$}

For $p = 7$, we derive a family of genus $3$ hyperelliptic 
curves
$$
	X_7: v^2 = u(u^3-1)^2 + 3t,
$$
and an endomorphism $\eta$ of $\Jac(X_7)$ 
with $\ZZ[\eta] \cong \ZZ[\zeta_7 + \zeta_7^{-1}]$ 
by Proposition~\ref{dihedralendomorphisms}.
Applying Algorithm~\ref{TNalgo},
we derive polynomial maps $T$ and $N$,
which we use with Algorithm~\ref{evaluationalgo}
to realize $\eta$ as $\eta: P \mapsto \textsc{Evaluate}(P,T,N)$.

\section{Applications II: Curves with Cyclotomic Covering}
\label{cyclotomicsection}

In this section we develop explicit endomorphisms for the one 
dimensional families of hyperelliptic curves with real multiplication 
based on cyclotomic coverings, as defined in Tautz, Top, and 
Verberkmoes~\cite{TTV-1991}.

\subsection{Construction of the Cyclotomic Covering}

Let $n \ge 2$, and let $\rho_n$ and $\rho_{2n}$ be primitive
$n^\mathrm{th}$ and $2n^\mathrm{th}$ roots of unity over $k$ such that
$\rho_{2n}^2 = \rho_n$; also set $\tau_n = \rho_n + \rho_n^{-1}$.
Consider the family of hyperelliptic curves of genus $n$ over $k$ 
in one free parameter $t$ defined by
$$
	C_n: y^2 = x(x^{2n}+tx^n + 1).
$$
The curve $C_n$ has 
an automorphism $\zeta$ of order $2n$
and an involution $\sigma$, defined by
$$
	\zeta: (x,y) \longmapsto (\rho_nx,\rho_{2n}y)
	\quad \mbox{and}\quad 
	\sigma: (x,y) \longmapsto \left(x^{-1}, x^{-(n+1)}y\right),
$$
respectively;
note that $\zeta^n$ is the hyperelliptic involution 
$(x,y)\mapsto (x,-y)$.  
We define $X_n := C_n/\langle\sigma\rangle$ to be the quotient of 
$C_n$ by the action of $\sigma$.
The curve $X_n$ has an an affine model
$$
	X_n: v^2 = f_n(u) = D_n(u,1) + t, 
$$
where $D_n(u,1)$ is the $n^\mathrm{th}$ Dickson polynomial 
of the first kind with parameter\footnote{
	Dickson polynomials are generally defined with a parameter $a$
	in $k$, by the recurrence
	$$
		D_n(u,a) = uD_{n-1}(u,a) - aD_{n-2}(u,a).
	$$
	It is easily shown that the curve defined by 
	$v^2 = D_n(u,a) + t$ for any nonzero $a$ is a twist of $X_n$.
	When $a = 0$,
	we obtain a one-dimensional family of curves
	with complex multiplication by $\ZZ[\zeta_n]$;
	these curves are described in \cite[\S6.4]{Smith-thesis}.
}~$1$,
defined recursively by
\begin{equation}
\label{dicksondefinition}
	D_n(u,1) = uD_{n-1}(u,1) - D_{n-2}(u,1)
\end{equation}
for $n \ge 2$,
with $D_0(u,1) = 2$ and $D_1(u,1) = u$.
Dickson polynomials and their properties are described in
\cite{LidlMullenTurnwald-1993};
for our purposes, 
it is enough to know that
\begin{equation}
\label{dicksonidentity1}
	D_n(u + u^{-1}, 1) = u^n + u^{-n}
\end{equation}
(this is easily verified by induction),
which further implies
\begin{equation}
\label{dicksonidentity2}
	D_{nm}(u,1) = D_n(D_m(u,1),1) .
\end{equation}

\begin{remark}
	When $n$ is odd,
	our curves $C_n$ and $X_n$ coincide with the
	curves $\mathcal{D}_n$ and $\mathcal{C}_n$ of \cite{TTV-1991};
	for even $n$,
	our families instead coincide 
	with the curves described in
	the remark of \cite[page 1058]{TTV-1991}.
\end{remark}

\noindent
The quotient projection $\pi: C_n \to X_n$
is a covering of degree $2$. 
Equation~\eqref{dicksonidentity1} above
shows that it is defined by 
$$
	\pi: (x,y) \longmapsto (u,v) = (x+x^{-1}, x^{-(n+1/2)}y).
$$
The automorphism $\zeta$ of $C_n$
induces an endomorphism $\eta = (\pi\circ\zeta)_*\circ\pi^*$ of
$\Jac(X_n)$. If $n$ is prime,
then Proposition~\ref{dihedralendomorphisms} 
implies that $\ZZ[\eta] \cong \ZZ[\zeta_n + \zeta_n^{-1}]$,
where $\zeta_n$ is an $n^\mathrm{th}$ root of unity over $\QQ$.

The endomorphism $\eta$ 
is induced by the correspondence 
$Z := (\pi\circ\zeta\times\pi)(C_n)$
on $X_n\times X_n$,
for which we directly compute an affine model
$$
	Z = \variety{
		v_2 - v_1,
		u_2^2 + u_1^2 - \tau_n u_1 u_2 + \tau_n^2 - 4
	}.
$$
Setting $t_1 := \tau_n u$
and $n_1 := u^2 + \tau_n^2 - 4$,
we apply Algorithm~\ref{TNalgo} to obtain maps
$T: k[u] \to k[u]$ and $N: k[u] \to k[u]$
such that the endomorphism $\eta$ is realized by
$P \mapsto \textsc{Evaluate}(P,T,N)$,
using Algorithm~\ref{evaluationalgo}.
The first few $t_i$ and $n_{i,j}$ derived in Algorithm~\ref{TNalgo}
are given in Table~\ref{cyclotomictsandns} below.

\begin{table}
\caption{Cyclotomic covers: $t_i$ and $n_{i,j}$ for $0 \le i \le j \le
3$}\vspace{-4mm}
\label{cyclotomictsandns}
\begin{center}
\begin{tabular}{|c|@{\;}l|}
\hline
$t_0$	& $2$ 		\\
$t_1$	& $\tau_n u$	\\
$t_2$	& $(\tau_n^2 - 2)u^2 - 2(\tau_n^2 - 4)$ \\
$t_3$ 	& $\tau_n(\tau_n^2 - 3)u^3 - 3\tau_n(\tau_n^2 - 4)u$ \\
\hline
$n_{0,0}$ & $1$ \\
$n_{0,1}$ & $\tau_n u$ \\
$n_{0,2}$ & $(\tau_n^2 - 2)u^2 - 2(\tau_n^2 - 4)$ \\
\hline
\end{tabular}
\begin{tabular}{|c|@{\;}l|}
\hline
$n_{0,3}$ & $\tau_n((\tau_n^2 - 3) u^2 - 3(\tau_n^2 - 4))u$\\
$n_{1,1}$ & $u^2 + \tau_n^2 - 4$ \\
$n_{1,2}$ & $\tau_n (u^2 + \tau_n^2 - 4) u$ \\
$n_{1,3}$ & $(\tau^2 - 2)u^4 + (\tau^2 - 4)^2(u^2 - 2)$\\
$n_{2,2}$ & $(u^2 + \tau_n^2 - 4)^2$\\
$n_{2,3}$ & $\tau_n (u^2 + \tau_n^2 - 4)^2) u$ \\
$n_{3,3}$ & $(u^2 + \tau_n^2 - 4)^3$\\
\hline
\end{tabular}
\end{center}
\vspace{-4mm}
\end{table}

The elliptic curve $C_1: y^2 = x(x^2+tx+1)$ 
is obviously covered by $C_n$,
and is therefore a factor of $\Jac(C_n)$.
The following analogue of 
Theorem~\ref{prop_endo_artin-schreier} holds for this cyclotomic 
family, and is proved similarly.

\begin{proposition}
\label{prop_endo_cyclo}
The Jacobian $\Jac(C_n)$ is isogenous to $C_1 \times \Jac(X_n)^2$ 
for $n$ prime, and its endomorphism ring contains an order in 
$\QQ \times \MM_2(\QQ(\eta_n))$.
\end{proposition}

\begin{remark}
	If
	$n$ is a prime other than $5$,
	then \cite[Corollary 6]{TTV-1991} implies that 
	$\Jac(X_n)$ is absolutely simple for general values of $t$
	over a field of characteristic $0$.
	For $n = 5$,
	we find that the condition of Stoll \cite{Stoll-1995}
	(see \cite[\S14.4]{Flynn})
	for $\Jac(X_5)$ to be absolutely simple 
	is satisfied by $X_5$ with $t = 1$ at $p = 11$.
	Conversely, 
	if $n = pm$, for $p > 2$ and $m > 1$, 
	then identity~\eqref{dicksonidentity2} above
	gives a covering $X_n \to X_p$ of degree $m$,
	defined by $(u,v) \mapsto (D_m(u,1), v)$.
	It follows that $\Jac(X_n)$ 
	has a factor isogenous to $\Jac(X_p)$,
	and so is not simple.
\end{remark}

\subsection{Hyperelliptic Curves of Genus $2$ with Real Multiplication
by $\eta_5$}

Consider the case $n = 5$.  
Equation~\eqref{dicksondefinition}
shows that $D_5(u,1) = u^5 - 5u^3 + 5u$,
so the curve $X_5 = C_5/\subgrp{\sigma}$ 
is the curve of genus $2$ defined by the affine model
$$
	X_5: v^2 = f_5(u) = u^5 - 5u^3 + 5u + t.
$$

Each point on $\mathrm{Jac}(X_5)$ has a representative in the form
$(a(u),v-b(u))$, with $\deg a = 2$ and $\deg b = 1$;
so suppose $a(u) = a_2u^2 + a_1u + a_0$ and $b(u) = b_1u + b_0$.
Applying Algorithm~\ref{TNalgo},
we obtain maps $T$ and $N$ such that
$$
\begin{array}{rl}
	N(a) = & a_2^2u^4
		+ a_2a_1\tau_5u^3 
		+ (2a_2^2(\tau_5^2 - 4) + a_1^2 + a_2a_0(\tau_5^2 - 2))u^2 
	\\ &{}
		+ a_1(a_2(\tau_5^2 - 4) + a_0)\tau_5u 
		+ ((\tau_5^2 - 4)(a_2^2(\tau_5^2 - 4) + a_1^2 - 2a_2a_0) +
a_0^2),	\\
	N(b) = & b_1^2u^2 
		+ b_1b_0\tau_5u
		+ b_1^2(\tau_5^2 - 4) + b_0^2,  
	\mbox{ and }
	\\
	T(b) = & \tau_5 b_1 u + 2 b_0.
\end{array}
$$
The endomorphism $\eta$ is then explicitly realized by
$\eta: P \mapsto \textsc{Evaluate}(P,T,N)$,
using Algorithm~\ref{evaluationalgo}.

\begin{remark}
The weighted projective Igusa invariants of the generic curve are:
$$
(140:550:640t^2-60:22400t^2-77725:256t^4-2048t^2+4096).
$$
In particular,
this family corresponds to a one-dimensional subvariety 
in the moduli space.
\end{remark}

\subsection{Hyperelliptic Curves of Genus $3$ with Real Multiplication
by $\eta_7$}

%
In the case $n = 7$,
we derive a family of curves
$$
	X_7: v^2 = u^7 - 7u^5 + 14u^3 - 7u + t ,
$$
and an endomorphism $\eta$ of $\Jac(X_7)$
with $\ZZ[\eta] \cong \ZZ[\zeta_7 + \zeta_7^{-1}]$
by Proposition~\ref{dihedralendomorphisms}.
Applying Algorithm~\ref{TNalgo},
we derive polynomial maps $T$ and $N$,
which we may then use with Algorithm~\ref{evaluationalgo}
to realize $\eta$ as $\eta: P \mapsto \textsc{Evaluate}(P,T,N)$.


\section{Applications III: Curves From Elliptic Coverings}
\label{mestresection}

In \cite{Mestre-1991}, Mestre constructs a series of two dimensional 
families of hyperelliptic curves with explicit real endomorphisms, 
which are similarly realized by explicit correspondences. 
For the case $\eta_5$, Takashima~\cite{Takashima-2006} independently 
developed an explicit algorithm and complexity analysis for two and 
three dimensional families\footnote{%
  The moduli of genus 2 curves with real multiplication by $\eta_5$ 
  form a two dimensional subvariety of the moduli space of genus 2 
  curves, so this three dimensional family contains one dimensional 
  fibres of geometrically isomorphic curves.
  } 
referred to as Mestre--Hashimoto and Brumer--Hashimoto (see~\cite{Hash2000}).

\subsection{Hyperelliptic Curves of Genus $2$ with Real Multiplication
by $\eta_5$}

Let $s$ and $t$ be free parameters,
and consider the family of curves defined by
$$
	X_5 : v^2 = f_5(u) 
	= u^4(u-s) - s(u+1)(u-s)^3 + s^3u^3 - tu^2(u-s)^2.
$$
Mestre shows that $\Jac(X_5)$ has an endomorphism $\eta$ 
satisfying $\eta^2 + \eta - 1 = 0$,
induced by the correspondence $Z$
with affine model
$$
	Z = \variety{
		v_2 - v_1,
		u_1^2u_2^2 + s(s - 1)u_1u_2 - s^2(u_1 - u_2) + s^3
	}.
$$
We will derive
an explicit form for $\eta$.
Since $X_5$ is a curve of genus $2$,
each point of $\Jac(X_5)$
may be represented by an ideal $(a(u),v-b(u))$ 
with $a = a_2u^2+a_1u+a_0$
and $b = b_1u+b_0$.  
Setting $t_1 = -s((s - 1)u_2 - s)/u_2^2$
and $n_1 = s^2(u_2 + s)/u_2^2$,
we apply Algorithm~\ref{TNalgo}
to derive maps $T$ and $N$
such that
$$
\begin{array}{r@{\;=\;}l}
	N(a) & a_2^2n_{2,2} + a_2a_1n_{1,2} + a_1^2n_{1,1}
     		+ a_2a_0a_2n_{0,2} + a_1a_0n_{0,1} + a_0^2n_{0,0}, \\
	N(b) & b_1^2 n_{1,1} + b_1b_0 n_{0,1} + b_0^2 n_{0,0}, 
	\mbox{ and } \\
	T(b) & -b_1s((s - 1)u - s)/u^2 + 2b_0,
\end{array}
$$
with the $n_{i,j}$ given in the table below.
\begin{center}
\begin{tabular}{|c|@{\;}l|}
\hline
	$n_{0,0}$ & $1$\\
	$n_{0,1}$ & $-s((s - 1)u - s)/u^2$\\
	$n_{0,2}$ & $s^2(((s-1)u-s)^2 - 2u^2(u+s))/u^4$\\
\hline
\end{tabular}
\begin{tabular}{|c|@{\;}l|}
\hline
	$n_{1,1}$ & $s^2(u + s)/u^2 $\\
	$n_{1,2}$ & $-s^3(u + s)((s - 1)u - s)/u^4 $\\
	$n_{2,2}$ & $s^4(u+s)^2/u^4 $\\
\hline
\end{tabular}
\end{center}
The endomorphism $\eta$ is then explicitly realized by
$\eta: P \mapsto \textsc{Evaluate}(P,T,N)$,
using Algorithm~\ref{evaluationalgo}.

\subsection{Hyperelliptic Curves of Genus $3$ with Real Multiplication
by $\eta_7$}

Let $s$ and $t$ be free parameters, 
and consider the family 
of hyperelliptic genus $3$ curves 
defined by 
$$
	X_7 : v^2 = f_7(u) = \phi_7(u) - t\,\psi_7(u)^2
$$
where $\psi_7(u) := u(u-s^3+s^2)(u-s^2+s)$ and
$$
\begin{array}{c@{\,}l}
\phi_7(u) := u\psi_7(u)^2
    & +\, s(s-1)(s^2-s+1)(s^3+2s^2-5s+1)u^5\\
    & -\, s^3(s-1)^2(6s^4-11s^3+12s^2-11s-1)u^4\\
    & +\, s^4(s-1)^3(s^2-s-1)(s^3+2s^2+6s+1)u^3\\
    & -\, s^6(s-1)^4(s+1)(3s^2-5s-3)u^2\\
    & +\, s^8(s-1)^5(s^2-3s-3)u + s^{10}(s-1)^6.
\end{array}
$$
Mestre shows that $\Jac(X_7)$ 
has an endomorphism $\eta$ 
satisfying 
$\eta^3 + \eta^2 - 2\eta - 1 = 0$,
induced by the correspondence
$ Z = \variety{ v_2 - v_1, E }$
on $X_7\times X_7$,
where
$$
	E = 
	u_1^2u_2^2 - s^2(s-1)(s^2-s-1)u_1u_2 - s^4(s-1)^2(u_1+u_2) 
		+ s^6(s-1)^3.
$$
Since $X_7$ is a curve of genus $3$,
each point on $\Jac(X_7)$ may be represented by an ideal
$(a(u),v-b(u))$,
where $a$ and $b$ are polynomials of degree $3$ and $2$, 
respectively. 
Setting
$$\begin{array}{r@{\;=\;}l}
	t_1 	& s^2(s-1)((s^2-s-1)u + s^2(s-1)) / u^2 
	\quad \mbox{and} \\
	n_1 	& -s^4(s-1)^2(u + s^2(s-1)) / u^2 ,
\end{array} $$
we apply Algorithm~\ref{TNalgo} 
to derive maps $T$ and $N$ from $k[u]$ into $k(u)$;
the elements $n_{i,j}$ computed by Algorithm~\ref{TNalgo}
are given in the table below.
\begin{center}
\begin{tabular}{|c|l|}
\hline
	$n_{0,0}$ & $1$\\
	$n_{0,1}$ & $s^2(s-1)((s^2-s-1)u + s^2(s-1))/u^2$\\
	$n_{0,2}$ & $(s - 1)^2s^4(2u^3 + (s^4 - 3s^2 + 2s + 1)u^2$\\
          	& $\quad +\, 2(s-1)(s^2 - s - 1)s^2u +
(s-1)^2s^4)/u^4$\\	$n_{1,1}$ & $-s^4(s-1)^2(u + s^2(s-1))/u^2$\\
	$n_{1,2}$ & 
 	$(s^6(s-1)^3(s^2 - s - 1)u^2 + s^9(s - 1)^5u +
s^{10}(s-1)^5)/u^4$\\	$n_{2,2}$ & $s^8(s-1)^4(u + s^2(s-1))^2/u^4$\\
\hline
\end{tabular}
\end{center}
The endomorphism $\eta$ is then explicitly realized by
$\eta: P \mapsto \textsc{Evaluate}(P,T,N)$,
using Algorithm~\ref{evaluationalgo}.

\section{Construction of Curves of Cryptographic Proportions} 

The curves presented here not only admit efficiently computable 
endomorphisms, they also permit random selection of curve 
parameters in a large family.  For example, let $\FF_{5^{37}} 
= \FF_5[\xi]$ be extension of $\FF_5$ such that $\xi^{37} + 
4\xi^2 + 3\xi + 3 = 0$, and take
$$
	t = 3\xi^5 + \xi^4 + 3\xi^3 + \xi^2 + 2\xi + 3.
$$
This gives a curve $X: v^2 = u(u^2-1)^2+t$ in the Artin--Schreier 
family whose Jacobian has nearly prime group order
$$
	|\Jac(X)(\FF_5[\xi])| = 5\cdot n,
$$
with prime cofactor 
$$
n = 1058791184067701689674637025340531565456011790341311.
$$ 
Such curves are amenable to efficient point counting 
techniques using Monsky-Washnitzer
cohomology~\cite{GG-2003,Kedlaya-2001}.  If $y$ is a square root of $t$,
then $(0,y)$ is a point on $X$; let $P = [(u,v-y)]$ be the corresponding
point on $J$.  Then $Q = [5](P)$ generates a cyclic group of order $n$, 
on which $[\eta]$ satisfies 
$$
	([\eta_5]^2 + [\eta_5]-1)(Q) = [(1)]
$$
and in particular, $[\eta_5](5P) = [m](5P)$, where
$$
	m = 336894053941004885519266617028956898972619907667301
$$
is one of the two roots of $x^2+x-1 \bmod n$.

\vspace{4mm}

\noindent
{\bf Acknowledgement.}  The authors thank K.~Takashima for providing 
an advance draft of his article~\cite{Takashima-2006}, 
and for references 
to the work of Hashimoto.

\end{document}